\documentclass[12pt]{article}
\usepackage{amssymb, epsfig,amssymb, latexsym}
\usepackage{amsfonts,psfrag,amsmath,bbm,color,multirow,algorithm,url}
\usepackage{placeins}
\usepackage{amscd}

\graphicspath{{./figs/}}

\newtheorem{remark}{Remark}[section]

\newtheorem{definition}{Definition}[section]

\def\PP{{{\rm l}\kern - .15em {\rm P} }}
\def\PN2{{\PP_{N}-\PP_{N-2}}}

\newcommand{\R}{\mathbbm{R}}

\newcommand{\cO}{\mathcal{O}}

\newcommand{\bphi}{\boldsymbol{\varphi}}

\newcommand{\btau}{\boldsymbol{\tau}}

\newcommand{\ba}{\boldsymbol{a}}
\newcommand{\bas}{{\boldsymbol a}^{snap}}

\newcommand{\bff}{\boldsymbol{f}}
\newcommand{\bFF}{{\boldsymbol F}}
\newcommand{\bG}{{\boldsymbol G}}
\newcommand{\bGs}{{\boldsymbol G}^{snap}}

\newcommand{\bu}{\boldsymbol{u}}

\newcommand{\bur}{{\boldsymbol{u}}_r}

\newcommand{\bx}{\boldsymbol{x}}
\newcommand{\bX}{\boldsymbol{X}}

\newcommand{\bXr}{{\bf X}^r}

\newcommand{\tA}{\tilde{A}}
\newcommand{\tB}{\widetilde{B}}

\newcommand{\obu}{\overline{\boldsymbol u}}

\newcommand{\as}{a^{snap}}

\newcommand{\Gs}{G^{snap}}

\newcommand{\deleted}[1]{{}}%\color{red}{#1}}}
\usepackage{fullpage}

\begin{document}

\title{Calibrated Filtered Reduced Order Modeling}

\author{
X. Xie
\thanks{Department of Mathematics, Virginia Tech, Blacksburg, VA 24061
Partially supported by NSF DMS1522656, email: xupingxy@vt.edu}
\and M. Mohebujjaman
\thanks{Department of Mathematical Sciences, Clemson University, Clemson, SC, 29634;
Partially supported by NSF DMS1522191, email: mmohebu@g.clemson.edu}
\and L. G. Rebholz
\thanks{Department of Mathematical Sciences, Clemson University, Clemson, SC 29634; 
Partially supported by NSF DMS1522191, email: rebholz@clemson.edu}
\and T. Iliescu
\thanks{Department of Mathematics, Virginia Tech, Blacksburg, VA 24061
Partially supported by NSF DMS1522656, email: iliescu@vt.edu}.
}

\date{\today}

\maketitle

\begin{abstract}
We propose a calibrated filtered reduced order model (CF-ROM) framework for the numerical simulation of general nonlinear PDEs that are amenable to reduced order modeling.
The novel CF-ROM framework consists of two steps:
(i) In the first step, we use explicit ROM spatial filtering of the nonlinear PDE to construct a filtered ROM.
This filtered ROM is low-dimensional, but is not closed (because of the nonlinearity in the given PDE).
(ii) In the second step, we use a calibration procedure to close the filtered ROM, i.e., to model the interaction between the resolved and unresolved modes.
To this end, we use a linear or quadratic ansatz to model this interaction and close the filtered ROM.
To find the new coefficients in the closed filtered ROM, we solve an optimization problem that minimizes the difference between the full order model data and our ansatz.  
Although we use a fluid dynamics setting to illustrate how to construct and use the CF-ROM framework, we emphasize that it is built on general ideas of spatial filtering and optimization and is  independent of (restrictive) phenomenological arguments. 
Thus, the CF-ROM framework can be applied to a wide variety of PDEs.
\end{abstract}



\section{Introduction}
    \label{sec:introduction}

Reduced order models (ROMs) have been successfully used to reduce the computational cost of scientific and engineering applications that are governed by relatively few recurrent dominant spatial structures~\cite{ballarin2016fast,bistrian2015improved,gunzburger2017ensemble,hesthaven2015certified,HLB96,noack2011reduced,perotto2017higamod,quarteroni2015reduced,stefanescu2015pod}.
For a given general partial differential equation (PDE), the ROM strategy for approximating its solution ($\bu$), is straightforward:
(i) choose modes $\{ \bphi_1, \ldots, \bphi_d \}$, which represent the recurrent spatial structures of the given PDE;
(ii) choose the dominant modes $\{ \bphi_1, \ldots, \bphi_r \}$, $r \leq d$, as  basis functions for the ROM;
(iii) use a Galerkin projection of the given PDE onto the ROM space $\bXr := \text{span} \{ \bphi_1, \ldots, \bphi_r \}$ to obtain a low-dimensional dynamical system, which represents the ROM;  
(iv) in an offline stage, compute the ROM vectors, matrices and tensors; and 
(v) in an online stage, repeatedly use the ROM (for various parameter settings and/or longer time intervals).

The quest for better ROMs has yielded two types of approaches: 
In the first approach, the ROM is fixed and the ROM basis is improved (see, e.g., the recent survey in~\cite{taira2017modal}).  
In the second approach, the basis is fixed and the ROM is improved.
In this paper, we exclusively use the latter.
Specifically, we try to answer the following question:
\begin{equation}
	\boxed{\text{For a fixed ROM space $\bXr$,  what is the best ROM in $\bXr$?}}
	\label{eqn:introduction-1}
\end{equation}

\vspace*{0.1cm}

The standard ROM is developed as follows:
(i) the given nonlinear PDE($\bu$) is projected onto the ROM space $\bXr$;
(ii) the Galerkin truncation $\bur = \sum_{j=1}^{r} u_j \, \bphi_j$ is used;
(iii) $\bu$ is replaced with $\bur$ in step (i).
The standard ROM is often efficient and relatively accurate~\cite{ballarin2015supremizer,HLB96,noack2011reduced}, but can fail in realistic applications
when large numbers of modes are needed to accurately represent the system. 
To ensure a low computational cost, it is desirable that ROMs use only a few modes (i.e., low $r$ values) and discard the remaining modes $\{ \bphi_{r+1}, \ldots, \bphi_d \}$.
The resulting low-dimensional ROM, however, can yield inaccurate results and numerical oscillations (see, e.g.,~\cite{amsallem2012stabilization,balajewicz2013low,barone2009stable,benosman2016robust,bergmann2009enablers,carlberg2017galerkin,giere2015supg,osth2014need,protas2015optimal}).
The general explanation for these inaccurate results is that the models fail to account for the interaction between resolved and unresolved modes~\cite{AHLS88,CSB03,gouasmi2016characterizing,noack2008finite,wang2012proper,wells2017evolve,xie2017approximate}.

To address the standard ROM's inaccuracy, we propose a novel {\it calibrated-filtered ROM (CF-ROM)} framework for the numerical simulation of {\it general nonlinear PDEs}.
The new framework explicitly models the interaction between the resolved  modes $\{ \bphi_1, \ldots, \bphi_r \}$ and unresolved modes $\{ \bphi_{r+1}, \ldots, \bphi_d \}$.
This is in stark contrast with the standard ROM, which simply ignores this interaction. 
We believe that the new CF-ROM framework represents a potential answer to the question~\eqref{eqn:introduction-1}, as we explain below, after first describing the general ideas.
The CF-ROM framework is schematically illustrated in~\eqref{eqn:cf-rom-schematic}, below.

\begin{equation}
\boxed{
\begin{CD}
	\text{FOM}		
	@>
	\text{filtering} 
	>>	
	\text{F-ROM}
	@>
	\text{calibration} 
	>>	
	\text{CF-ROM} \\[0.1cm]
\end{CD}
}
	\label{eqn:cf-rom-schematic}
\end{equation}
In the first step of \eqref{eqn:cf-rom-schematic}, the full order model (FOM), which contains all the information in the underlying system, is filtered with an {\it explicit ROM spatial filter}.
The resulting {\it filtered ROM (F-ROM)} approximates only the large spatial structures of the system and, therefore, requires fewer modes than the FOM.
To create a usable model, one still needs to solve the {\it ROM closure problem} to model the interaction between the resolved and unresolved modes in the F-ROM, which is done in the second step of \eqref{eqn:cf-rom-schematic}.  
We use a calibration procedure to {\it construct} the actual ROM (i.e., not to data fit its solution~\cite{abou2016new,cordier2010calibration,couplet2005calibrated,wang20162d}). 
We employ a linear or quadratic ansatz to model the interaction between the resolved and unresolved modes in the F-ROM, and then we find the ROM coefficients by solving a (low-dimensional) optimization problem that minimizes the difference between the full order model data and our ansatz.

Our proposed CF-ROM framework has some connections to some other popular models, in particular the {\it nonlinear Galerkin}~\cite{FMRT01} and {\it variational multiscale}~\cite{hughes1998variational} methods, since they all use the small-large scale separation. 
The CF-ROM is also reminiscent of large eddy simulation (LES) in a sense, since it uses spatial filtering.
However, CF-ROM is different from the nonlinear Galerkin method since it can use general types of spatial filtering, whereas the latter is based on projection. 
CF-ROM is also different from the variational multiscale method and LES since it uses calibration to model the interaction with the unresolved modes, whereas the other two methods use phenomenological arguments based on, e.g., the energy cascade and eddy viscosity.  
Indeed, the fact that CF-ROM does not use phenomenological arguments is precisely whey it is not restricted in its application, and can be utilized for a wide variety of nonlinear PDEs.

\vspace*{0.3cm}

The rest of the paper is organized as follows:
In Section~\ref{sec:rom}, we present the POD and the standard ROM.
In Section~\ref{sec:f-rom}, we introduce the filtered ROM.
In Section~\ref{sec:c-rom}, we present the general ROM calibration procedure.
In Section~\ref{sec:cf-rom}, we use calibration to solve the closure problem in the filtered ROM and to construct the CF-ROM.
In Section~\ref{sec:numerical-results}, we present a numerical investigation of the CF-ROM in the simulation of the one-dimensional Burgers equation.  
Finally, in Section~\ref{sec:conclusions}, we draw conclusions and outline future research directions.  

\section{Reduced Order Models (ROMs)}
    \label{sec:rom}

In this section, we briefly review POD (Section~\ref{sec:pod}) and the standard ROM (Section~\ref{sec:g-rom}).  Although the new CF-ROM framework can be applied to many types of nonlinear PDEs, to present the method it is necessary to pick a particular model, and so we select our favorite, the incompressible Navier-Stokes equations (NSE):
\begin{eqnarray}
    && \frac{\partial \bu}{\partial t}
    - Re^{-1} \Delta \bu
    + \bu \cdot \nabla \bu
    + \nabla p
    = {\bf 0} \, ,
    \label{eqn:nse-1}                                                         \\
    && \nabla \cdot \bu
    = 0 \, ,
    \label{eqn:nse-2}
\end{eqnarray}
where $\bu$ is the velocity, $p$ the pressure and $Re$ the Reynolds number. 
We use the initial condition $\bu(\bx, 0) = \bu_0(\bx)$ and (for simplicity) homogeneous Dirichlet boundary conditions: $\bu(\bx, t) = {\bf 0}$.

    \subsection{Proper Orthogonal Decomposition (POD)}
        \label{sec:pod}
        One of the most popular reduced order modeling techniques is the
        POD~\cite{HLB96,noack2011reduced,Sir87abc}. For the snapshots
        $\{\bu^1_h,\ldots, \bu^{N_s}_h\}$, which are, e.g., finite element (FE) solutions of
        \eqref{eqn:nse-1}--\eqref{eqn:nse-2} at $N_s$ different time instances,
        the POD seeks a low-dimensional basis that approximates the snapshots
        optimally with respect to a certain norm. In this paper, the commonly
        used $L^2$-norm will be chosen. The solution of the minimization
        problem is equivalent to the solution of the eigenvalue problem
        $
            YY^TM_h \bphi_j = \lambda_j \bphi_j,
            \  j=1,\ldots,N,
        $
        where $\bphi_j$ and $\lambda_j$ denote the vector of the FE coefficients
        of the POD basis functions and the POD eigenvalues, respectively, $Y$
        denotes the snapshot matrix, whose columns correspond to the FE
        coefficients of the snapshots, $M_h$ denotes the FE mass matrix, and $N$
        is the dimension of the FE space $\bX^h$.
        The eigenvalues are real and non-negative, so they can be ordered as
        follows: $\lambda_1 \ge \lambda_2 \ge \ldots \ge \lambda_d \ge \lambda_{d + 1} = \ldots = \lambda_N = 0$, where $d$ is the rank of the snapshot matrix.
        The POD basis consists of the normalized functions $\{
        \bphi_{j}\}_{j=1}^{r}$, which correspond to the first $r\le N$ largest
        eigenvalues. Thus, the POD space is defined as $\bX^r := \text{span} \{
        \bphi_1, \ldots, \bphi_r \}$.

\subsection{Standard Galerkin ROM (G-ROM)}
    \label{sec:g-rom}

The POD approximation of the velocity is defined as 
\begin{equation}
	{\bu}_r({\bf x},t) 
	\equiv \sum_{j=1}^r a_j(t) \bphi_j({\bf x}) \, ,
	\label{eqn:g-rom-1}
\end{equation}
where $\{a_{j}(t)\}_{j=1}^{r}$ are the sought time-varying coefficients, which are determined by solving the {\it Galerkin ROM (G-ROM)}: $\forall \, i = 1, \ldots, r,$
    \begin{eqnarray}
        \left(
            \frac{\partial \bu_r}{\partial t} , \bphi_{i}
        \right)
        + Re^{-1} \, \left(
            \nabla \bu_r ,
            \nabla \bphi_{i}
        \right)
        + \biggl(
            (\bu_r \cdot \nabla) \bu_r ,
            \bphi_{i}
        \biggr)
        = 0 \, .
    \label{eqn:g-rom}
    \end{eqnarray}
In~\eqref{eqn:g-rom}, we assumed that the modes $\{ \bphi_1, \ldots, \bphi_r \}$ are perpendicular to the discrete pressure space, which is the case if standard, conforming LBB stable elements (such as Taylor-Hood, Scott-Vogelius, the mini-element, etc.) are used for the snapshot creation.
Plugging~\eqref{eqn:g-rom-1} into~\eqref{eqn:g-rom} yields the following low-dimensional dynamical system:  $\forall \, i = 1 \ldots r$
\begin{equation}
	\dot{a}_i	
	 = \sum_{m=1}^{r} A_{i m} \, a_m(t)
	 + \sum_{m=1}^{r} \sum_{n=1}^{r} B_{i m n} \, a_n(t) \, a_m(t) \, ,
	\label{eqn:g-rom-2}
\end{equation}
where
\begin{eqnarray}
	A_{im}
	= - Re^{-1} \, \left( \nabla \bphi_m , \nabla \bphi_i \right),
	\qquad
	B_{imn}
	= - \bigl( \bphi_m \cdot \nabla \bphi_n , \bphi_i \bigr) \, .
	\label{eqn:g-rom-3b}
\end{eqnarray}

\section{Filtered ROMs (F-ROMs) and the ROM Closure Problem}
	\label{sec:f-rom}
	
In this section, we present details regarding the first step of the CF-ROM framework~\eqref{eqn:cf-rom-schematic}.  
That is, we discuss the creation of the {\it filtered ROM (F-ROM)}, from which the ROM closure problem will reveal itself.  
To this end, in Section~\ref{sec:explicit-rom-spatial-filtering} we present the ROM projection, which is the explicit ROM spatial filter that we use to construct the CF-ROM.
In Section~\ref{sec:f-rom-framework}, we develop the F-ROM framework, which has also been used to develop LES-ROMs~\cite{xie2017approximate}.
Finally, in Section~\ref{sec:rom-closure-problem}, we outline the celebrated ROM closure problem, which needs to be solved in the F-ROM.
We emphasize that the ROM closure problem is treated completely differently in CF-ROM and LES-ROM: CF-ROM uses calibration, while LES-ROMs use phenomenological arguments (e.g., energy cascade and eddy viscosity). 

\subsection{Explicit ROM Spatial Filtering}
    \label{sec:explicit-rom-spatial-filtering}

Spatial filtering has been used in ROMs, mainly as a preprocessing tool to eliminate the noise in the snapshot data (see, e.g., Section 5 in~\cite{aradag2011filtered} for a survey of relevant work).  
However, our approach is fundamentally different:
We explicitly use spatial filtering in the construction of the actual ROM, not in the development of the ROM basis.  
In this paper, we exclusively use the ROM projection~\cite{wang2012proper,wells2017evolve} as  spatial filter, but we note that we could also other spatial filters (e.g., the ROM differential filter~\cite{sabetghadam2012alpha,wells2017evolve,xie2017approximate}).

For a fixed $r \leq d$ and a given $\bu \in \bX^h$, the ROM projection seeks $\obu^r \in \bX^{r}$ such that
        \begin{eqnarray}
            \left( \obu^r, \bphi_j \right)
            = ( \bu , \bphi_j )
            \quad \forall \, j=1, \ldots r \, .
            \label{eqn:rom-projection}
        \end{eqnarray}
The ROM projection~\eqref{eqn:rom-projection} has been used for theoretical purposes, e.g., in the error analysis of the G-ROM~\cite{KV01}, and as a smoothing ROM spatial filter in~\cite{wang2012proper,wells2017evolve}.

\subsection{F-ROM Framework}
\label{sec:f-rom-framework}

We now outline the F-ROM framework, which has also been used to  develop LES-ROMs in~\cite{xie2017approximate}.
To this end, we use the standard LES approach~\cite{layton2012approximate,rebollo2014mathematical,Sag06}, which consists of the following steps:
(i) Use an explicit spatial filter to filter the NSE;
(ii) Use the resulting spatially filtered NSE and the ROM approximation to obtain the F-ROM.

Using the fact that $\nabla \cdot \bu = 0$ in~\eqref{eqn:nse-1}, we can write the nonlinear term as follows:  $( \bu \cdot \nabla) \, \bu = \nabla \cdot (\bu \, \bu)$.
Filtering the NSE, assuming that differentiation and filtering commute~\cite{Sag06}, and projecting the equations onto a space of weakly divergence-free functions $\boldsymbol{\phi}$, we obtain the spatially filtered NSE (see equations (35)--(36) in~\cite{xie2017approximate}):
\begin{eqnarray}
	\left(
            \frac{\partial \obu}{\partial t} , \boldsymbol{\phi}
        \right)
        + Re^{-1} \, \biggl(
            \nabla \obu ,
            \nabla \boldsymbol{\phi}
        \biggr)
        + \biggl(
            \nabla \cdot \bigl( \obu \, \obu \bigr) ,
            \boldsymbol{\phi}
        \biggr)
        + \biggl(
            \nabla \cdot \btau ,
            \boldsymbol{\phi}
        \biggr)
        = {\bf 0} \, ,
    \label{eqn:les-rom-1}
\end{eqnarray}
where
\begin{eqnarray}
\btau
= {\overline{\bu \, \bu}} - \obu \, \obu
\label{eqn:les-rom-2}
\end{eqnarray}
is the {\it subfilter-scale stress tensor}.
    
The spatial structures in the spatially filtered NSE~\eqref{eqn:les-rom-1} are larger than the spatial structures in the NSE~\eqref{eqn:nse-1}.
Thus, we expect that for a fixed target numerical accuracy of the ROM, the spatially filtered NSE will require fewer POD modes than the NSE, which is advantageous from a computational point of view.

Of course, to develop a useful ROM from the spatially filtered NSE~\eqref{eqn:les-rom-1}, we need to do the following:
(i) use a ROM approximation for the continuous velocity field $\obu$; and
(ii) use a ROM approximation of the spatial filter.
We make the following choices:
Use $\bu_d \sim \bu$ in (i) (where $d$ is the rank of the snapshot matrix) and $\obu^r \sim \obu$ in (ii) (where $\obu^r$ is the ROM projection~\eqref{eqn:rom-projection}).

Using $\bu_d \sim \bu$ in (i) means that we employ the best possible approximation of the continuous velocity field $\bu$ in the set of snapshots (i.e., in $\bX^d$). 
Using the ROM projection onto $\bX^r$ as the spatial filter in (ii) means that we are projecting the equations from $\bX^d$ onto $\bX^r$.
With these choices in the spatially filtered NSE~\eqref{eqn:les-rom-1}, we obtain the spatially {\it filtered ROM (F-ROM)}:
$\forall \, i = 1, \ldots, r,$
\begin{eqnarray}
	\left(
            \frac{\partial \overline{\bu_d}^r}{\partial t} , \bphi_i
        \right)
        + Re^{-1} \, \biggl(
            \nabla  \overline{\bu_d}^r ,
            \nabla \bphi_i
        \biggr)
        + \biggl(
            \nabla \cdot \bigl( \overline{\bu_d}^r \, \overline{\bu_d}^r \bigr) ,
            \bphi_i
        \biggr)
        + \biggl(
            \nabla \cdot \btau_r ,
            \bphi_i
        \biggr)
        = {\bf 0} \, ,
    \label{eqn:f-rom-0}
\end{eqnarray}
where
\begin{eqnarray}
\btau_r
= \overline{\bu_d \, \bu_d}^r - \overline{\bu_d}^r \, \overline{\bu_d}^r \, .
\label{eqn:les-rom-4}
\end{eqnarray}
Since we are using the ROM projection onto $\bX^r$ as the spatial filter, we have
\begin{eqnarray}
\overline{\bu_d}^r
= \bu_r .
\label{eqn:les-rom-4b}
\end{eqnarray}
Thus, the F-ROM~\eqref{eqn:f-rom-0} becomes
\begin{eqnarray}
	\left(
            \frac{\partial \bu_r}{\partial t} , \bphi_i
        \right)
        + Re^{-1} \, \biggl(
            \nabla \bu_r ,
            \nabla \bphi_i
        \biggr)
        + \biggl(
            \nabla \cdot \bigl( \bu_r \, \bu_r \bigr) ,
            \bphi_i
        \biggr)
        + \biggl(
            \nabla \cdot \btau_r ,
            \bphi_i
        \biggr)
        = {\bf 0} \, ,
    \label{eqn:f-rom}
\end{eqnarray}
where the ROM stress tensor is
\begin{eqnarray}
\btau_r
= \overline{\bu_d \, \bu_d}^r - \bu_r \, \bu_r \, .
\label{eqn:les-rom-6}
\end{eqnarray}

\subsection{ROM Closure Problem}
    \label{sec:rom-closure-problem}

The F-ROM~\eqref{eqn:f-rom} is an $r$-dimensional ODE system for $\bur$.
Since $r \ll N$, the F-ROM~\eqref{eqn:f-rom} is a computationally efficient surrogate model for the FOM (i.e., the FE approximation of the NSE, which is an $N$-dimensional ODE system).
We emphasize, however, that the F-ROM~\eqref{eqn:f-rom} is not a closed system of equations, since the ROM stress tensor $\btau_r$ depends on $\bu_d$ (see equation~\eqref{eqn:les-rom-6}).
Thus, to close the F-ROM~\eqref{eqn:f-rom}, we need to solve the ROM closure problem~\cite{galletti2007accurate,noack2008finite,sirisup2004spectral,wang2012proper}.

\begin{definition}[ROM closure problem]
The ROM closure problem is modeling the ROM stress tensor $\btau_r$, i.e., finding a function $\bff$ that satisfies
\begin{equation}
	\nabla \cdot \btau_r
	\approx \bff(\ba , \bX^r) \, .
	\label{eqn:les-rom-7}
\end{equation}
The ROM closure model is the explicit dependence of the ROM stress tensor $\btau_r$ on the ROM coefficients $\ba$ and ROM space $\bX^r$ in equation~\eqref{eqn:les-rom-7}.
	\label{def:rom-closure-problem}
\end{definition}
Note that by neglecting the last term on the LHS of~\eqref{eqn:f-rom}, the F-ROM~\eqref{eqn:f-rom} is identical to the standard G-ROM~\eqref{eqn:g-rom}.
Thus, {\it formally}, one could write the following decomposition:
\begin{equation}
	\boxed{\text{F-ROM = G-ROM + ROM Closure Model}}
	\label{eqn:les-rom-9}
\end{equation}

\vspace*{0.1cm}

In the ROM community, it is well known that the standard G-ROM with few ROM modes (i.e., with small $r$) generally yields inaccurate results in complex applications.
Probably the most popular ``solution" to this problem is to supplement the G-ROM with extra terms:
\begin{equation}
	\boxed{\text{Good ROM = G-ROM + Extra Terms}}
	\label{eqn:f-rom-10}
\end{equation}

\vspace*{0.1cm}

Numerical and physical arguments are used to build the Extra Terms in~\eqref{eqn:f-rom-10}.
Furthermore, the decomposition~\eqref{eqn:f-rom-10} is usually called ``solving the ROM closure problem"~\cite{galletti2007accurate,noack2008finite,sirisup2004spectral,wang2012proper},
and the Extra Terms are thought to model the effect of the discarded ROM modes $\{ \bphi_{r+1}, \ldots, \bphi_{d} \}$.
We emphasize, however, that to our knowledge no precise definition of the ROM closure problem has been given so far.
The F-ROM~\eqref{eqn:f-rom} allows for the first time a clear definition of the ROM closure problem (in Definition~\ref{def:rom-closure-problem}).

\begin{remark}[F-ROM Consistency]
	Probably the most popular ROMs of the form~\eqref{eqn:f-rom-10} are the eddy viscosity closure models pioneered by Lumley and his group~\cite{AHLS88,HLB96} and extended to state-of-the-art LES closure models in~\cite{wang2012proper}.
	As emphasized in, e.g., Section 2.4 in~\cite{balajewicz2016minimal}, although these eddy viscosity ROM closure models have been successfully applied to a large number of complex, high Reynolds number flows, they have a significant drawback, namely the loss of consistency between the NSE and the ROM. 
	Indeed, since the G-ROM in~\eqref{eqn:f-rom-10} is modified empirically, the resulting eddy viscosity ROM no longer corresponds to a Galerkin projection of the NSE. 
	
	The F-ROM, however, is consistent with the NSE.
	Indeed, the F-ROM~\eqref{eqn:f-rom} is just the projection of the full order model (i.e., the best representation of the NSE in the snapshot space, $\bX^d = \text{span} \{ \bphi_1, \ldots, \bphi_d \}$) onto the ROM space, $\bX^r$.
	\label{remark:f-rom-consistency}
\end{remark}

\section{Calibrated ROMs (C-ROMs)}
    \label{sec:c-rom}

In this section, we outline the general ROM calibration procedure (see, e.g.,~\cite{cordier2010calibration}), and in Section~\ref{sec:cf-rom} we use it to construct the new CF-ROM from the F-ROM.

\subsection{Modified G-ROM}
	\label{sec:modified-g-rom}

We start by rewriting the standard G-ROM~\eqref{eqn:g-rom-2} as 
\begin{equation}
	\dot{\ba}
	 = \bFF(\ba) \, ,
	\label{eqn:c-rom-1}
\end{equation}
where, for $i = 1, \ldots, r$,
\begin{equation}
	 F_{i}(\ba)
	 = b_i
	 + \sum_{m=1}^{r} A_{i m} \, a_m
	 + \sum_{m=1}^{r} \sum_{n=1}^{r} B_{i m n} \, a_n \, a_m \, .
	\label{eqn:c-rom-2}
\end{equation}
In the ROM literature, ROM calibration is used when, for different reasons, the G-ROM~\eqref{eqn:g-rom-2} is deemed inaccurate and is replaced by the following {\it modified G-ROM}:
\begin{equation}
	\dot{\ba}
	 = \bFF(\ba)
	 + \bG \, .
	\label{eqn:c-rom-3}
\end{equation}
Thus, the modified G-ROM~\eqref{eqn:c-rom-3} depends on a new quantity, $\bG$, whose aim is to make the modified G-ROM~\eqref{eqn:c-rom-3} more accurate than the standard G-ROM~\eqref{eqn:g-rom-2}.
Examples of quantities $\bG$ used in the current C-ROMs include the pressure representation proposed in~\cite{noack2005need}.

\subsection{Ansatz: Linear or Quadratic?}
	\label{sec:ansatz}

To make the modified G-ROM~\eqref{eqn:c-rom-3} fit the general ROM framework, the following ansatz is generally used in the current C-ROMs:
\begin{equation}
	\bG
	= \bG(\ba) \, .
	\label{eqn:c-rom-4}
\end{equation}
In fact, to make the modified G-ROM~\eqref{eqn:c-rom-3} even further resemble the standard G-ROM~\eqref{eqn:g-rom-2}, the following choices are generally made in ansatz~\eqref{eqn:c-rom-4}: For $i = 1, \ldots, r,$
\begin{eqnarray}
	&& \hspace*{-1.5cm}  G_i(\ba)
	= \sum_{m=1}^{r} \tA_{i m} \, a_m
	\hspace*{4.40cm} \text{linear ansatz}  
	\label{eqn:c-rom-5} \\[0.3cm]
	&& \hspace*{-1.5cm}  \text {or} \nonumber \\[0.3cm]
     && \hspace*{-1.5cm} G_i(\ba)
	= \sum_{m=1}^{r} \tA_{i m} \, a_m 
	+ \sum_{m=1}^{r} \sum_{n=1}^{r} \tB_{i m n} \, a_m \, a_n
	\hspace*{0.5cm} \text{quadratic ansatz} \, . 
	\label{eqn:c-rom-6}
\end{eqnarray}
Note that the linear ansatz~\eqref{eqn:c-rom-5} is a particular case of the quadratic ansatz~\eqref{eqn:c-rom-6}, in which $\tB_{i m n} = 0, \ \forall \, m , n$.  
Thus, for clarity purposes, in this section we will exclusively use the quadratic ansatz~\eqref{eqn:c-rom-6} as our base ansatz.

The linear ansatz~\eqref{eqn:c-rom-5} was used, e.g., in~\cite{buffoni2006low} (for the ROM closure and ROM pressure, see page 145), \cite{galletti2007accurate} (for the ROM closure and ROM pressure, see pages 356--357),  \cite{galletti2004low} (for the ROM pressure, see page 164) and~\cite{noack2005need} (for the ROM pressure, see page 361). 
The quadratic ansatz~\eqref{eqn:c-rom-6} was used, e.g., in~\cite{noack2005need} (for the ROM pressure, see equations (2.18)--(2.19)).
The linear ansatz and the quadratic ansatz were shown to yield similar qualitative results for the pressure representation (see page 361 in~\cite{noack2005need}).

Using the quadratic ansatz~\eqref{eqn:c-rom-6} in the modified G-ROM~\eqref{eqn:c-rom-3} yields, for $i = 1, \ldots, r,$
\begin{equation}
	\dot{a_i}
	 = b_i
	 + \sum_{m=1}^{r} \left( A_{i m} + \tA_{i m} \right) \, a_m
	 + \sum_{m=1}^{r} \sum_{n=1}^{r} \left( B_{i m n} + \tB_{i m n} \right) \, a_n \, a_m \, .
	\label{eqn:c-rom-7}
\end{equation}

\subsection{Optimization: Constrained or Unconstrained?}
	\label{sec:optimization}

The ROM calibration aims at finding the coefficients $\tA_{i m}$ and $\tB_{i m n}$ that ensure the highest accuracy of the modified G-ROM~\eqref{eqn:c-rom-7}.
In the current C-ROMs, the accuracy is based on the difference between the ROM coefficients $a_i$ in the modified G-ROM~\eqref{eqn:c-rom-7} and the ROM coefficients obtained from the snapshots, $\as_i$.
The values $\as_i(t_j)$, with snapshot time instances $t_j , \ j = 1, \ldots, M$, are obtained by projecting the corresponding snapshots $\bu(t_j) = \sum_{i=1}^{d} \as_i(t_j) \, \bphi_i$ onto the POD basis functions $\bphi_i$ and using the orthogonality of the POD basis functions: For $i = 1, \ldots, d$ and $j = 1, \ldots, M,$
\begin{equation}
	\as_i(t_j)
	= \biggl( \bu(t_j) , \bphi_i \biggr) \, .
	\label{eqn:c-rom-8}
\end{equation}

To find the coefficients $\tA_{i m}$ and $\tB_{i m n}$ that ensure the highest accuracy of the modified G-ROM~\eqref{eqn:c-rom-7}, the current C-ROMs solve an optimization problem, which can be of two types: an ODE {\it constrained optimization} problem
\begin{eqnarray}
	&& \min_{\tA_{i m}, \, \tB_{i m n}} \sum_{j = 1}^{M} \| \ba(t_j) - \bas(t_j) \|^2
	\label{eqn:c-rom-9} \\
	&& \text{s.t. } \ba \text { solves } \eqref{eqn:c-rom-7} \, ,
	\label{eqn:c-rom-10}
\end{eqnarray}
where $\| \cdot \|$ is the Euclidian norm in $\R^r$, or an {\it unconstrained optimization} problem
\begin{eqnarray}
	\min_{\tA_{i m}, \, \tB_{i m n}} \sum_{j = 1}^{M} \| \bG(\bas(t_j)) - \bGs \|^2
	\label{eqn:c-rom-11} \, ,
\end{eqnarray}
where $\bGs$ is the {\it true} quantity $\bG$ computed with snapshot data.

The constrained minimization problem~\eqref{eqn:c-rom-9}--\eqref{eqn:c-rom-10} was used, e.g., in~\cite{buffoni2006low} (see page 145), \cite{cordier2010calibration} (see page 277), \cite{couplet2005calibrated} (see page 203), \cite{galletti2007accurate} (see page 357), and~\cite{galletti2004low} (see page 164).
As explained on page 357 in~\cite{galletti2007accurate}, the constrained minimization problem~\eqref{eqn:c-rom-9}--\eqref{eqn:c-rom-10} is usually solved in two steps: 
First, one passes from the discrete to the continuous time variable by using a spline interpolation.
Then,  in the continuous setting, one solves an optimality problem that consists of the direct problem, the adjoint problem and optimality conditions.

The unconstrained minimization problem~\eqref{eqn:c-rom-11} was used, e.g., in~\cite{noack2005need}  (see page 361).
This unconstrained minimization problem was solved in~\cite{noack2005need} at a discrete level (not at a continuous level, as done for the constrained minimization problem~\eqref{eqn:c-rom-9}--\eqref{eqn:c-rom-10}).
Linear regression was used in~\cite{noack2005need} to solve the resulting overdetermined linear system.

\paragraph{Comparison}

Both the constrained minimization problem~\eqref{eqn:c-rom-9}--\eqref{eqn:c-rom-10} and the unconstrained minimization problem~\eqref{eqn:c-rom-11} have been used to determine the C-ROM coefficients.
For example, Noack et al.~\cite{noack2005need} solved the unconstrained minimization problem~\eqref{eqn:c-rom-11}, whereas Galletti et al.~\cite{galletti2004low} solved the constrained minimization problem~\eqref{eqn:c-rom-9}--\eqref{eqn:c-rom-10}.

We are not aware of any comparison of the two approaches.
Although a comparison of the two approaches is needed, it is beyond the scope of the present paper.
We only mention that, for a fair comparison, we believe that the {\it predictive capabilities} of the C-ROM should be assessed (see, e.g., \cite{wang2012proper} for other examples of predictive ROMs).
Indeed, we believe that by solving the constrained minimization problem~\eqref{eqn:c-rom-9}--\eqref{eqn:c-rom-10} (as in Galletti et al.~\cite{galletti2004low}), one effectively ensures that the C-ROM results are as close as possible to the NSE results.
We emphasize, however, that this is true only on the time interval on which the snapshots were collected, say, $[0,T]$.
On a longer time interval, say, $[0,10 \,T]$, the C-ROM obtained with the approach used in~\cite{galletti2004low} might not perform optimally.
In other words, the predictive capabilities of the C-ROM are not clear. 
On the other hand, by solving the unconstrained minimization problem~\eqref{eqn:c-rom-11} (as in Noack et al.~\cite{noack2005need} and in this paper), one aims at finding the best model for the unknown quantity.
The hope is that this model has some universal features that can be determined from the available snapshots.
If this is the case, then the C-ROM obtained with the approach used in Noack et al.~\cite{noack2005need} and in this paper might have predictive capabilities, i.e., it might yield accurate results on a longer time interval, say, $[0,10 \,T]$.
The following quote from page 361 in Noack et al.~\cite{noack2005need} seems to support our conjecture:
``In a recent study of wake flows, Galletti et al.~\cite{galletti2004low} use the same ansatz for the pressure term. 
In that study the coefficients are determined from an optimization problem which minimizes the error between the Galerkin solution and the Navier--Stokes simulation. 
The present empirical approach minimizes the error between the pressure term in primitive variables and in the Galerkin model representation. 
By construction, the term of Galletti et al.~\cite{galletti2004low} is more accurate in the short term whereas the present term is more forgiving of long-term phase errors."
	
For simplicity purposes, in this paper we are using the unconstrained minimization problem~\eqref{eqn:c-rom-11}. 	
We plan to compare the two approaches and investigate the potential predictive capabilities of the CF-ROM in a future study.

\subsection{C-ROM Algorithm}
	\label{sec:c-rom-algorithm}

The optimal coefficients $\tA^{opt}_{i m}$ and $\tB^{opt}_{i m n}$ that solve the minimization problems~\eqref{eqn:c-rom-9}--\eqref{eqn:c-rom-10} or~\eqref{eqn:c-rom-11} are used in the modified G-ROM~\eqref{eqn:c-rom-7}, yielding the {\it calibrated ROM (C-ROM)}
\begin{equation}
	\dot{a_i}
	 = \sum_{m=1}^{r} \left( A_{i m} + \tA^{opt}_{i m} \right) \, a_m
	 + \sum_{m=1}^{r} \sum_{n=1}^{r} \left( B_{i m n} + \tB^{opt}_{i m n} \right) \, a_n \, a_m \, .
	\label{eqn:c-rom-12}
\end{equation}

In Algorithm~\ref{alg:rom-calibration}, we summarize the C-ROM for a generic quantity of interest, $\bG$.

\begin{algorithm}[H]
	\caption{C-ROM}
	\label{alg:rom-calibration}
  	\begin{enumerate}
		\item[Step 1] Consider the modified G-ROM $\dot{\ba} = \bFF(\ba) + \bG$
    		\item[Step 2] Use ansatz $\bG = \bG(\ba)$ for the unknown quantity $\bG$:
\begin{eqnarray}
	&& \hspace*{-1.5cm}  G_i(\ba)
	= \sum_{m=1}^{r} \tA_{i m} \, a_m
	\hspace*{4.4cm} \text{linear ansatz}  
	\label{eqn:c-rom-5b} \\[0.3cm]
	&& \hspace*{-1.5cm}  \text {or} \nonumber \\[0.3cm]
     && \hspace*{-1.5cm} G_i(\ba)
	= \sum_{m=1}^{r} \tA_{i m} \, a_m 
	+ \sum_{m=1}^{r} \sum_{n=1}^{r} \tB_{i m n} \, a_m \, a_n
	\hspace*{0.5cm} \text{quadratic ansatz} \, . 
	\label{eqn:c-rom-6b}
\end{eqnarray}
		\item[Step 3] Solve the constrained minimization problem~\eqref{eqn:c-rom-9b}--\eqref{eqn:c-rom-10b} or the unconstrained minimization problem~\eqref{eqn:c-rom-11b} to determine the unknown coefficients $\tilde{A}_{i m}$ and $\tilde{B}_{i m n}$ in~\eqref{eqn:c-rom-5b}--\eqref{eqn:c-rom-6b}:
\begin{eqnarray}
	&& \min_{\tA_{i m}, \, \tB_{i m n}} \sum_{j = 1}^{M} \| \ba(t_j) - \bas(t_j) \|^2
	\label{eqn:c-rom-9b} \\
	&& \text{s.t. } \ba \text { solves } \eqref{eqn:c-rom-7} \, ,
	\label{eqn:c-rom-10b}
\end{eqnarray}
or
\begin{eqnarray}
	\min_{\tA_{i m}, \, \tB_{i m n}} \sum_{j = 1}^{M} \| \bG(\bas(t_j)) - \bGs \|^2
	\label{eqn:c-rom-11b} \, .
\end{eqnarray}
		\item[Step 4] The C-ROM has the following form:
		\begin{equation}
	\dot{a_i}
	 = \sum_{m=1}^{r} \left( A_{i m} + \tA^{opt}_{i m} \right) \, a_m
	 + \sum_{m=1}^{r} \sum_{n=1}^{r} \left( B_{i m n} + \tB^{opt}_{i m n} \right) \, a_n \, a_m \, ,
	\label{eqn:c-rom-12b}
\end{equation}
	where the optimal coefficients $\tA^{opt}_{i m}$ and $\tB^{opt}_{i m n}$ solve the minimization problems~\eqref{eqn:c-rom-9b}--\eqref{eqn:c-rom-10b} or~\eqref{eqn:c-rom-11b}.
	\end{enumerate}
\end{algorithm}

\section{Calibrated Filtered ROM (CF-ROM)}
    \label{sec:cf-rom}

The F-ROM~\eqref{eqn:f-rom} naturally fits into the general ROM calibration framework developed in Section~\ref{sec:c-rom}.
Indeed, the F-ROM~\eqref{eqn:f-rom} can be written as the weak form of the modified G-ROM~\eqref{eqn:c-rom-3}, i.e., as the sum of the standard G-ROM~\eqref{eqn:g-rom} and the unknown ROM stress tensor contribution $\bigl( \nabla \cdot \btau_r , \bphi_i \bigr)$, where $\btau_r$ is defined in~\eqref{eqn:les-rom-6}.
Thus, the general ROM calibration procedure presented in Section~\ref{sec:c-rom} can be used as ROM closure modeling, i.e., to model the term $\bigl( \nabla \cdot \btau_r , \bphi_i \bigr)$ in~\eqref{eqn:f-rom} as a function of the ROM coefficients ($\ba$).
In the remainder of this section, we use the general ROM calibration framework developed in Section~\ref{sec:c-rom} to solve the closure problem in the F-ROM~\eqref{eqn:f-rom}.

We start by denoting by $\bG$ the unknown ROM stress tensor contribution: For $i = 1, \ldots, r,$
\begin{equation}
	G_i(\ba)
	:= \bigl( \nabla \cdot \btau_r , \bphi_i \bigr) \, .
	\label{eqn:c-les-rom-1}
\end{equation} 

Next, we make a linear or quadratic ansatz for $\bG$ in~\eqref{eqn:c-les-rom-1} (see equations~\eqref{eqn:c-rom-5} and \eqref{eqn:c-rom-6}, respectively):
\begin{equation}
	G_i
	= \sum_{m=1}^{r} \tA_{i m} \, a_m 
	+ \sum_{m=1}^{r} \sum_{n=1}^{r} \tB_{i m n} \, a_m \, a_n \, . 
	\label{eqn:c-les-rom-2}
\end{equation}

The next step in the ROM calibration framework is solving the constrained minimization problem~\eqref{eqn:c-rom-9}--\eqref{eqn:c-rom-10} or the unconstrained minimization problem~\eqref{eqn:c-rom-11}.
For simplicity purposes, we will solve the latter.
To this end, we first compute the ROM coefficients from the snapshot data ($\as_i(t_j)$) by using equation~\eqref{eqn:c-rom-8}.
Next, we compute the {\it true} quantity $\bG$ from the snapshot data: $\forall \, i = 1, \ldots ,r,$
\begin{eqnarray}
	\Gs_i(t_j)
	&=& \left( \nabla \cdot \left[ \overline{\bu^{snap}_d(t_j) \, \bu^{snap}_d(t_j)}^r - \bu^{snap}_r(t_j) \, \bu^{snap}_r(t_j) \right] , \bphi_i \right)
	\nonumber \\
	&=& \biggl(  \nabla \cdot \biggl[ \overline{ \left( \sum_{k_1=1}^{d} a^{snap}_{k_1}(t_j) \bphi_{k_1} \right) \, \left( \sum_{k_2=1}^{d} a^{snap}_{k_2}(t_j) \bphi_{k_2} \right) }^r 
	\nonumber \\
	&& -  \left( \sum_{k_3=1}^{r} a^{snap}_{k_3}(t_j) \bphi_{k_3} \right)  \, \left( \sum_{k_4=1}^{r} a^{snap}_{k_4}(t_j) \bphi_{k_4} \right) \biggr] , \bphi_i \biggr) \, ,
	\label{eqn:c-les-rom-3}
\end{eqnarray}
where
\begin{eqnarray}
	\bu^{snap}_d(t_j)
	= \sum_{k=1}^{d} a^{snap}_{k}(t_j) \, \bphi_{k} \, ,
	\qquad
	\bu^{snap}_r(t_j)
	= \sum_{k=1}^{r} a^{snap}_{k}(t_j) \, \bphi_{k} \, .
	\label{eqn:c-les-rom-3c}
\end{eqnarray}

\begin{remark}[Computational Efficiency]
	\label{remark:c-les-rom-2}
	In practical settings, where the rank of the snapshot matrix can be extremely large (e.g., $d= \mathcal{O}(1000)$), using $\bu^{snap}_d \in \bX^d$ in~\eqref{eqn:c-les-rom-3} would be extremely costly.
	One possible solution would be to replace $\bu^{snap}_d$ with, say, $\bu^{snap}_{2 r}$ and the ROM projection on $\bX^d$ with the ROM projection on $\bX^{2 r}$.
	We numerically investigate the accuracy of this approximation in Section~\ref{sec:numerical-results}.
\end{remark}

Finally, we compute the quantity $\bG$ from the linear or quadratic ansatz~\eqref{eqn:c-les-rom-2} and the snapshot data:
\begin{equation}
	G_i\bigl(\bas(t_j)\bigr)
	= \sum_{m=1}^{r} \tA_{i m} \, a^{snap}_m(t) 
	+ \sum_{m=1}^{r} \sum_{n=1}^{r} \tB_{i m n} \, a^{snap}_m(t) \, a^{snap}_n(t) \, . 
	\label{eqn:c-les-rom-4}
\end{equation}

Plugging~\eqref{eqn:c-rom-8},~\eqref{eqn:c-les-rom-3}, and~\eqref{eqn:c-les-rom-4} into the unconstrained minimization problem~\eqref{eqn:c-rom-11}, we obtain a {\it linear least squares} problem (if we use the linear ansatz in~\eqref{eqn:c-les-rom-4}) or a {\it nonlinear least squares} problem (if we use the quadratic ansatz in~\eqref{eqn:c-les-rom-4}).

Using the quadratic ansatz in~\eqref{eqn:c-les-rom-4} seems more appropriate from a physical point of view.
Indeed, the ROM stress tensor $\btau_r$ is closely connected to the nonlinearity in the NSE.
Since this nonlinearity yields a quadratic term in the G-ROM~\eqref{eqn:g-rom-2}, one would expect that the ROM stress tensor be best approximated by a quadratic term, too, just as in the quadratic ansatz in~\eqref{eqn:c-les-rom-4}.
From a numerical point of view, however, the resulting nonlinear least squares problem seems more challenging.
In this case, one could use, e.g., Gauss-Newton methods or similar approaches (see Chapter 10 in~\cite{nocedal2006numerical}).

Using the linear ansatz in~\eqref{eqn:c-les-rom-4} does not have as much physical support as using the quadratic ansatz.
From a numerical point of view, however, solving the resulting linear least squares problem is straightforward (e.g., by using Matlab).
Thus, in what follows (including the numerical investigation in Section~\ref{sec:numerical-results}), we will use the linear ansatz.
We plan to investigate the more physical quadratic ansatz in a future study.

Next, we show how the unconstrained minimization problem can be numerically solved with the linear ansatz. 
Plugging~\eqref{eqn:c-les-rom-3} and \eqref{eqn:c-les-rom-4} with the linear ansatz into the unconstrained minimization problem~\eqref{eqn:c-rom-11}, we obtain the following minimization problem:
\begin{equation}
		\min_{\tA_{i m}} \, \sum_{j = 1}^{M} \, \sum_{i=1}^{r} \, \left| \, \sum_{m=1}^{r} \tA_{i m} \, a^{snap}_m(t_j) - \Gs_i(t_j) \, \right|^2 
	= \, \min_{\tA_{i m}} \, f(\tA_{i m}) \, .
	\label{eqn:c-les-rom-5}
\end{equation}
The optimality conditions for the minimization problem~\eqref{eqn:c-les-rom-5} yield
\begin{equation}
	\frac{\partial f}{\partial \tA_{k l}}
	= 0
	\qquad \forall \, k, l = 1, \ldots, r \, .
	\label{eqn:c-les-rom-6}
\end{equation}
The LHS of~\eqref{eqn:c-les-rom-6} can be calculated as follows:
\begin{eqnarray}
	\hspace*{-0.8cm}
	\frac{\partial f}{\partial \tA_{k l}}
	&=& \sum_{j = 1}^{M} 2 \, \left( \sum_{m=1}^{r} \, \tA_{k m} \, a^{snap}_m(t_j) - \bGs_k(t_j) \right) \, a^{snap}_l(t_j)
	\nonumber \\
	&=& 2 \, \sum_{j = 1}^{M} \left( \sum_{m=1}^{r} \, \tA_{k m} \, a^{snap}_m(t_j) \, a^{snap}_l(t_j) - \bGs_k(t_j) \, a^{snap}_l(t_j) \right) 
		\nonumber \\
	&=& \sum_{m=1}^{r} \, \tA_{k m} \, \sum_{j = 1}^{M} \biggl( a^{snap}_m(t_j) \, a^{snap}_l(t_j) \biggr)
	- \sum_{j = 1}^{M} \biggl( \bGs_k(t_j) \, a^{snap}_l(t_j) \biggr) .
	\label{eqn:c-les-rom-7}
\end{eqnarray}
Thus, the optimality conditions~\eqref{eqn:c-les-rom-6} reduce to the following matrix equation:
\begin{equation}
	\tA \, D
	= E \, ,
	\label{eqn:c-les-rom-8}
\end{equation}
where
\begin{eqnarray}
	D_{m l}
	:= \sum_{j = 1}^{M} \, a^{snap}_m(t_j) \, a^{snap}_l(t_j) \, ,
	\qquad
	E_{k l}
	:= \sum_{j = 1}^{M} \, \bGs_k(t_j) \, a^{snap}_l(t_j) \, .
	\label{eqn:c-les-rom-10}
\end{eqnarray}
Therefore, assuming that the matrix $D$ is invertible, the matrix $\tA$ is
\begin{equation}
	\tA 
	= E \, D^{-1} \, .
	\label{eqn:c-les-rom-11}
\end{equation}
The optimal coefficients $\tA^{CF-ROM}_{i m}$ that solve~\eqref{eqn:c-les-rom-11} yield the {\it calibrated  filtered ROM (CF-ROM)}:
\begin{eqnarray}
	\dot{a_i}
	 = \sum_{m=1}^{r} \left( A_{i m} + \tA^{CF-ROM}_{i m} \right) \, a_m
	 + \sum_{m=1}^{r} \sum_{n=1}^{r} B_{i m n} \, a_n \, a_m \, .
	\label{eqn:cf-rom}
\end{eqnarray}

The CF-ROM~\eqref{eqn:cf-rom} construction can be summarized in the following algorithm:

\begin{algorithm}[H]
	\caption{CF-ROM}
	\label{alg:cf-rom}
  	\begin{enumerate} \itemsep10pt
		\item Use~\eqref{eqn:c-rom-8} to compute $\as_i(t_j), \ i = 1, \ldots, d$.
		\item Use~\eqref{eqn:c-les-rom-10} to assemble matrix $D$.
    		\item Use~\eqref{eqn:c-les-rom-3} to compute $\Gs_i(t_j), \ i = 1, \ldots, r$.
			    Use the computationally efficient approach in Remark~\ref{remark:c-les-rom-2}.
		\item Use~\eqref{eqn:c-les-rom-10} to assemble matrix $E$.
    		\item Solve the matrix equation~\eqref{eqn:c-les-rom-11}.
    		\item The CF-ROM is given by~\eqref{eqn:cf-rom} with the optimal coefficients $\tA^{CF-ROM}_{i m}$ that solve the matrix equation~\eqref{eqn:c-les-rom-11} in step 5.
	\end{enumerate}
\end{algorithm}

\section{Numerical Results}
	\label{sec:numerical-results}

In this section, we numerically investigate the CF-ROM~\eqref{eqn:cf-rom} applied to the Burgers equation, which has become a common benchmark problem in ROM testing~\cite{borggaard2011artificial,KV01,san2014proper} due to its being a one-dimensional (1D) problem but yet displaying complex solution structure.  
 Indeed, Burgers equation is used to model phenomena in statistical mechanics, vehicle traffic, adhesion in cosmology, as well as simplistic modeling of fluid dynamics \cite{BK07}.
In Section~\ref{sec:burgers-equation}, we describe the mathematical and computational setting of the test problem.
In Section~\ref{sec:cf-rom-burgers}, we construct the CF-ROM for the Burgers equation, following a procedure that is analogous to that used above for the CF-ROM construction for the 3D NSE.
In Section~\ref{sec:cf-rom-vs-g-rom}, we compare the CF-ROM with the standard G-ROM, both in terms of numerical accuracy and computational efficiency.
Finally, in Section~\ref{sec:cf-rom-vs-state-of-the-art-roms}, we compare the CF-ROM with some of the most recent ROMs, both in terms of numerical accuracy and computational efficiency.

\subsection{Burgers Equation}
	\label{sec:burgers-equation}

For the test problem in the CF-ROM numerical investigation, we use the 1D Burgers equation as mathematical model, which was also used in~\cite{borggaard2011artificial,KV01,san2014proper}:
\begin{eqnarray}
  \left\{\begin{array}{ll}
  u_t-\nu \, u_{xx}+u \, u_x = 0& x\in [0,1],\\
  u(x,0)=u_0(x)& x\in [0,1],\\
  u(0,t)=u(1,t) = 0 \, , &  
  \end{array} \right.
\label{eqn:burgers}
 \end{eqnarray}
where $\nu=10^{-3}$ is the diffusion parameter and $t\in[0,1]$. 
We use the following initial condition:
\begin{eqnarray}
u_0(x)=\left\{\begin{array}{ll}
1&  x\in\left(0,\frac{1}{2}\right],\\
0&  x\in\left(\frac{1}{2},1\right).
\end{array}
\right.
\label{eqn:burgers-ic}
\end{eqnarray}

\begin{figure}[h]
\centering
\includegraphics[scale=0.5]{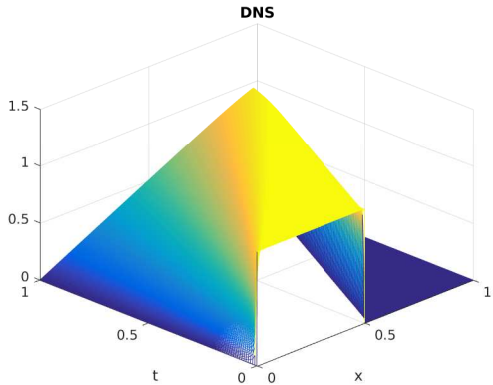}
\caption{Shown above is a plot of the solution of the Burgers equation DNS.}
\label{fig:burgers-dns}
\end{figure}

For the direct numerical simulation (DNS) (see Fig.~\ref{fig:burgers-dns}), we use a uniform mesh with $\Delta x=1/1024$ and a linear FE for the spatial discretization and the backward Euler method with a time step $\Delta t=10^{-3}$ for the time discretization. 
For all the ROMs, we use the forward Euler method with a time step $\Delta t=10^{-4}$ for the time discretization.
To generate the POD basis, we collect $101$ snapshots in the time interval $[0,1]$.

\subsection{CF-ROM for Burgers Equation}
	\label{sec:cf-rom-burgers}

In Section~\ref{sec:cf-rom}, the entire derivation of the CF-ROM is carried out for the general 3D NSE.
We now adapt the CF-ROM framework to the 1D Burgers equation~\eqref{eqn:burgers}, for which  the F-ROM~\eqref{eqn:f-rom} becomes the following:
$\forall \, i = 1, \ldots, r,$
\begin{eqnarray}
	\left(
            \frac{\partial u_r}{\partial t} , \varphi_i
        \right)
        + \nu \, \biggl(
            \frac{\partial u_r}{\partial x} ,
            \frac{\partial \varphi_i}{\partial x}
        \biggr)
        + \biggl(
            u_r \, \frac{\partial u_r}{\partial x} ,
            \varphi_i
        \biggr)
        + \biggl(
            f(\ba , X^r) ,
            \varphi_i
        \biggr)
        = 0 \, .
    \label{eqn:numerical-results-1}
\end{eqnarray}
Thus, $G_i(\ba)$ in~\eqref{eqn:c-les-rom-1} becomes the following:
$\forall \, i = 1, \ldots ,r,$
\begin{equation}
	G_i(\ba)
	:= \left( \, 
		\overline{ u_d \, \frac{\partial u_d}{\partial x} }^r 
		- u_r \, \frac{\partial u_r}{\partial x} \, , \,  
		\varphi_i \, 
	\right) \, .
	\label{eqn:numerical-results-2}
\end{equation} 
This, in turn, implies that~\eqref{eqn:c-les-rom-3} becomes:
$\forall \, i = 1, \ldots ,r,$
\begin{eqnarray}
	\Gs_i(t_j)
	&=& \biggl(  \biggl[ \overline{ \left( \sum_{k_1=1}^{d} a^{snap}_{k_1}(t_j) \varphi_{k_1} \right) \, \left( \sum_{k_2=1}^{d} a^{snap}_{k_2}(t_j) \frac{\partial \varphi_{k_2}}{\partial x} \right) }^r 
	\nonumber \\
	&& -  \left( \sum_{k_3=1}^{r} a^{snap}_{k_3}(t_j) \varphi_{k_3} \right)  \, \left( \sum_{k_4=1}^{r} a^{snap}_{k_4}(t_j) \frac{\partial \varphi_{k_4}}{\partial x} \right) \biggr] , \varphi_i \biggr) \, .
	\label{eqn:numerical-results-3}
\end{eqnarray}
To compute the first term in the inner product~\eqref{eqn:numerical-results-3}, we use the definition of the ROM projection~\eqref{eqn:rom-projection}:
$\forall \, j = 1, \ldots ,r,$
\begin{eqnarray}
	\left( 
		\sum_{k=1}^{r} b_k \, \varphi_k , \varphi_j 
	\right)
	= 
	\left( 
		\left( \sum_{k_1=1}^{d} a^{snap}_{k_1}(t_j) \varphi_{k_1} \right) \, \left( \sum_{k_2=1}^{d} a^{snap}_{k_2}(t_j) \frac{\partial \varphi_{k_2}}{\partial x} \right) , \varphi_j 
	\right) \, .
	\label{eqn:numerical-results-4}
\end{eqnarray}
Solving for $b_k, \, k = 1, \ldots, r$ in~\eqref{eqn:numerical-results-4}, we get 
\begin{eqnarray}
	\overline{ \left( \sum_{k_1=1}^{d} a^{snap}_{k_1}(t_j) \varphi_{k_1} \right) \, \left( \sum_{k_2=1}^{d} a^{snap}_{k_2}(t_j) \frac{\partial \varphi_{k_2}}{\partial x} \right) }^r
	=
	\sum_{k=1}^{r} b_k \, \varphi_k \, ,
	\label{eqn:numerical-results-5}
\end{eqnarray}
which allows us to compute $\Gs_i(t_j)$ in~\eqref{eqn:numerical-results-3}. 
At this stage, we have all the information needed to apply Algorithm~\ref{alg:cf-rom} and construct the CF-ROM~\eqref{eqn:cf-rom} for the Burgers equation.

\subsection{CF-ROM vs G-ROM}
	\label{sec:cf-rom-vs-g-rom}

\begin{figure}[H]
\centering
GROM (r=15) \hspace{2in} CF-ROM(r=15) \\
\includegraphics[width=0.45\textwidth,height=0.4\textwidth, trim=0 0 0 20, clip]{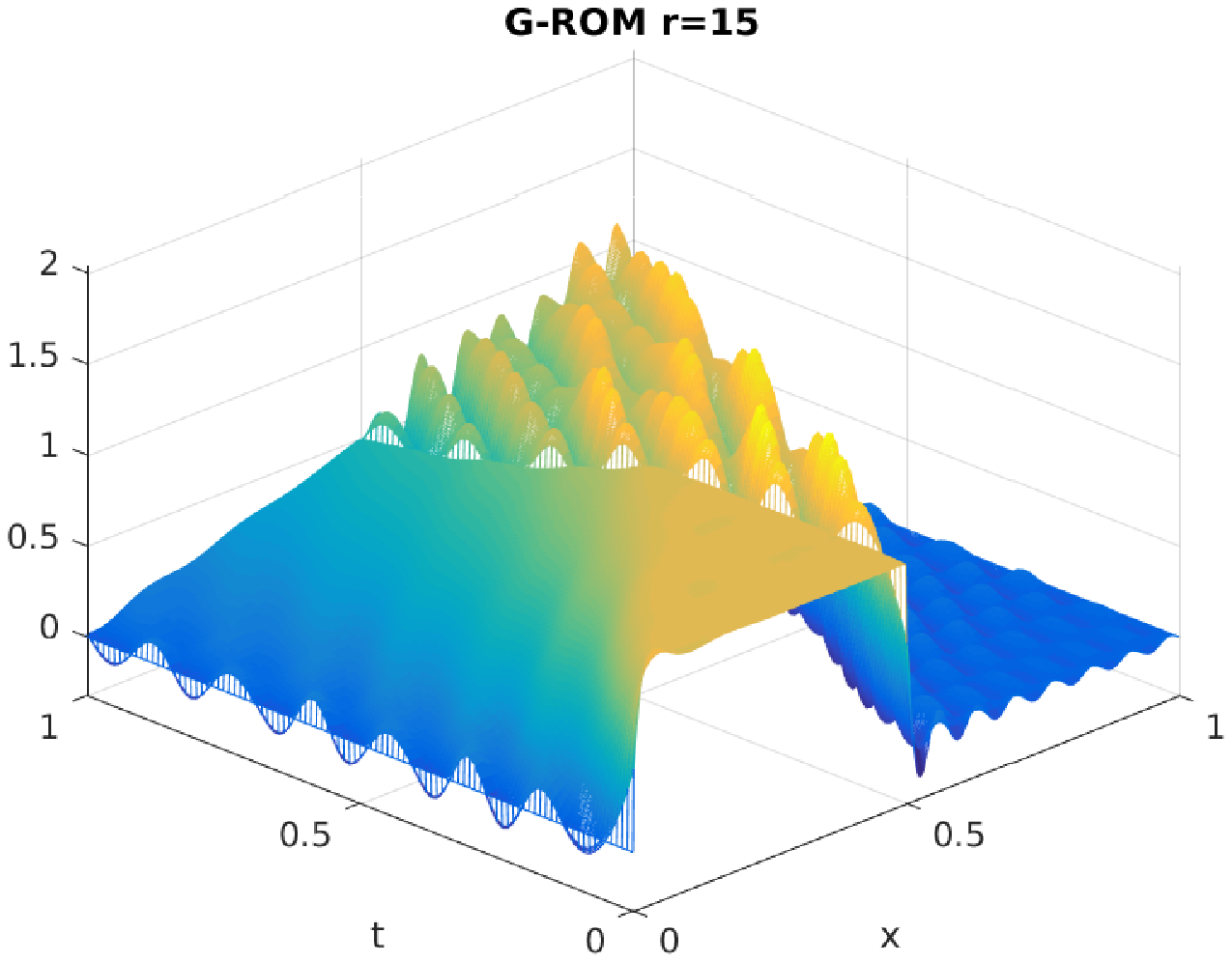}
\includegraphics[width=0.45\textwidth,height=0.4\textwidth, trim=0 0 0 20, clip]{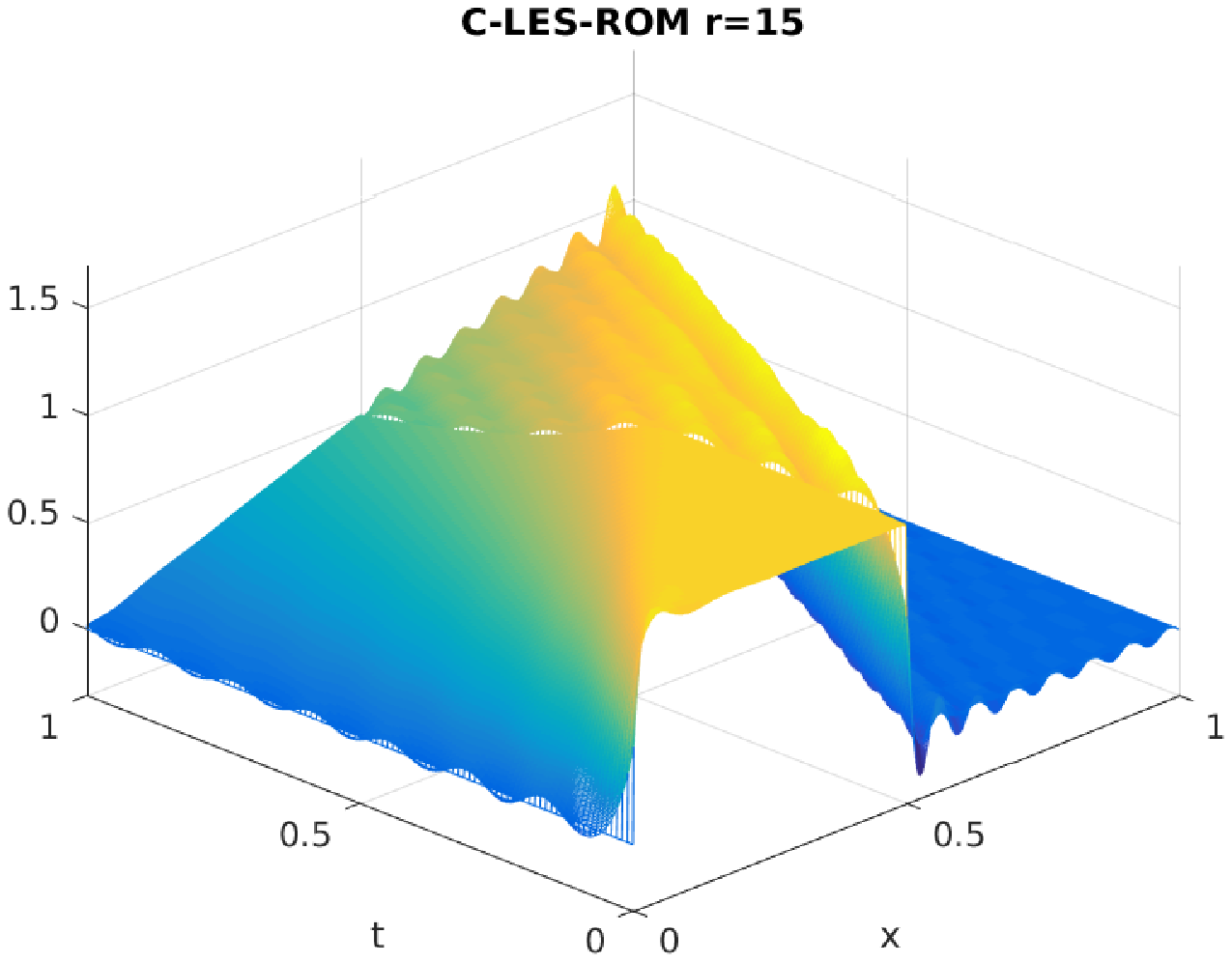}
\caption{
Shown above are plots for r=15 GROM (left) and CF-ROM (right) for the Burgers equation simulations.
	\label{fig:c-rom-vs-g-rom}
}
\end{figure}

For our first test, we compare the new CF-ROM~\eqref{eqn:cf-rom} with the standard G-ROM.
We consider the following $r$ values: $r=6, \, r=10$, and $r=15$.  
Table~\ref{table:c-rom-vs-g-rom} shows the maximum errors in these ROMs (comparing to the DNS solution in the $L^2$ norm at each time step, and then taking the average over the time steps).  
We observe that the CF-ROM error is consistently lower than the G-ROM error, generally by a factor of two.
To better illustrate the improvement of CF-ROM over G-ROM, we show in Figure~\ref{fig:c-rom-vs-g-rom} the respective solutions with $r=15$.

We also investigate the efficiency of methods and give timings for the various simulations in Table~\ref{table:c-rom-vs-g-rom}.
Although the CF-ROM is more accurate than the G-ROM, it is less efficient, and so a comparison is in order.  
In its original form, the CF-ROM uses~\eqref{eqn:c-les-rom-3} to compute $\Gs_i(t_j), \ i = 1, \ldots, r$.  
As mentioned in Remark~\ref{remark:c-les-rom-2}, this computation utilizes $\bu^{snap}_d \in \bX^d$.  
Since $d = \cO(1000)$ in some practical settings, this could make the CF-ROM computationally costly.  
Hence, following Remark~\ref{remark:c-les-rom-2}, we replace $\bu^{snap}_d$ in~\eqref{eqn:c-les-rom-3} with $\bu^{snap}_{m}$, where $r \leq m \leq d$:
\begin{equation}
	\bu^{snap}_d
	\approx
	\bu^{snap}_{m}
	\qquad r \leq m \leq d \, .
	\label{eqn:numerical-results-4}
\end{equation}
When $m$ in~\eqref{eqn:numerical-results-4} is low (i.e., close to $r$), the CF-ROM computational cost will be low, but the accuracy will also be low.
On the other hand, when $m$ in~\eqref{eqn:numerical-results-4} is large (i.e., close to $d$), the CF-ROM accuracy will be high, but the computational cost will also be high.  
Thus, we need to find an $m$ value in~\eqref{eqn:numerical-results-4} that ensures a compromise between accuracy and efficiency in the CF-ROM.  
From Table~\ref{table:c-rom-vs-g-rom}, we observe that $m=2r$ seems to achieve compromise between accuracy and efficiency in the CF-ROM.

\begin{table}[H]
\centering
\begin{tabular}{|l| l| c|c|c|c|c|c|c|}
\hline
Method &  Proj. Space & \multicolumn{3}{c|}{ $ \| u_{DNS} - u_{ROM} \| $ } &  \multicolumn{3}{c|}{CPU time (seconds)} \\ \hline
      &   & r=6 & r=10   & r=15     				&	r=6 & r=10   & r=15   \\ \hline
G-ROM &    & 0.2208& 0.1589&0.0848 			&	2.06 & 3.28 & 4.79 \\ \hline
CF-ROM & $X^r$   & 0.2208& 0.1589&0.0848 		&	2.06 & 3.28 & 4.79 \\ \hline
CF-ROM & $X^{r+1}$&  0.1214& 0.0854&0.0654 	&	2.16 & 3.43 & 4.95 \\ \hline
CF-ROM & $X^{2r}$ & 0.0928& 0.0627&0.0446 	&	2.75 & 4.42 & 6.67 \\ \hline
CF-ROM & $X^{3r}$ &  0.0923& 0.0632&0.0452 	&	3.44 & 5.61 & 8.57 \\ \hline
CF-ROM & $X^{d}$  &  0.0931& 0.0638&0.0454 	&	14.21 & 15.06 & 16.06 \\ \hline
\end{tabular}
\caption{
	G-ROM and CF-ROM errors and timings for $r=6, \, r=10$ and $r=15$ for Burgers equation simulation, using varying $r$ and ROM projection spaces.
	\label{table:c-rom-vs-g-rom}	}
\end{table}

\subsection{CF-ROM vs. State-Of-The-Art ROMs}
	\label{sec:cf-rom-vs-state-of-the-art-roms}
	
Above we have shown that the new CF-ROM is a clear improvement over the standard G-ROM, and a natural question is whether the CF-ROM is also an improvement over other, more accurate ROMs.
To address this question, we consider two recently proposed ROMs for fluid flows.
The first type of ROM that we consider is regularized ROMs (Reg-ROMs)~\cite{wells2017evolve}.
We consider two Reg-ROMs: the Leray ROM (L-ROM) and the evolve-then-filter ROM (EF-ROM).  
We consider a differential filter in both Reg-ROMs, since this ensured the most accurate results in~\cite{wells2017evolve}. 
The second type of ROM that we consider is an LES-ROM: the approximate deconvolution ROM (AD-ROM) proposed in~\cite{xie2017approximate}.

We compare the new CF-ROM with the AD-ROM, the L-ROM and the EF-ROM.
We test all ROMs on the Burgers equation test problem considered in this section.
We list the errors in Table~\ref{tab:rom-comparison-error}. 
The CF-ROM errors are slightly lower than all the other ROM errors.
We list the CPU times in Table~\ref{tab:rom-comparison-cpu-time}.
The CF-ROM CPU time is significantly lower than the CPU times of the other ROMs.
The results in Table~\ref{tab:rom-comparison-error} and Table~\ref{tab:rom-comparison-cpu-time} consistently show that, for this test problem, the CF-ROM is at least competitive with the other ROMs.

\begin{table}[h]
\centering
\begin{tabular}{|l|c|c|c|c|}
\hline
       & L-DF & EF-ROM & AD-ROM & CF-ROM($X^{2r}$) \\ \hline
$r=6$  & 0.1385 & 0.1005 & 0.1096 & 0.0928 \\ \hline
$r=10$ & 0.1135 & 0.0699 & 0.0633 & 0.0627 \\ \hline
$r=15$ & 0.1037 & 0.0549 & 0.0532 & 0.0446 \\ \hline
\end{tabular}
\caption{
	Burgers equation, ROM comparison:
	Errors for L-ROM-DF, EF-ROM-DF, AD-ROM and the new CF-ROM with $m=2r$ in~\eqref{eqn:numerical-results-4}.
	\label{tab:rom-comparison-error}
	}
\end{table}

\begin{table}[h]
\centering
\begin{tabular}{|l|c|c|c|c|}
\hline
       & L-DF & EF-ROM & AD-ROM & CF-ROM($X^{2r}$) \\ \hline
$r=6$  & 4.12 & 4.25 & 4.44 & 2.27 \\ \hline
$r=10$ & 6.72 & 6.91 & 7.26 & 4.42 \\ \hline
$r=15$ & 9.97 & 10.14 & 10.32 & 6.67 \\ \hline
\end{tabular}
\caption{
	Burgers equation, ROM comparison:
	CPU times for L-ROM-DF, EF-ROM-DF, AD-ROM and the new CF-ROM with $m=2r$ in~\eqref{eqn:numerical-results-4}.
	\label{tab:rom-comparison-cpu-time}
	}
\end{table}

\section{Conclusions and Outlook}
	\label{sec:conclusions}

In this paper, we proposed a novel ROM framework for the numerical simulation of general nonlinear PDEs.
This framework was based on explicit ROM spatial filtering and calibration.
The explicit ROM spatial filtering ensured computational efficiency of the filtered ROM.
The calibration procedure solved the ROM closure problem in the filtered ROM.
We numerically investigated the resulting CF-ROM in the simulation of the 1D Burgers equation.
First, we compared the new CF-ROM with the standard G-ROM.
The CF-ROM was significantly more accurate than the G-ROM.
Furthermore, the computational costs of the CF-ROM and G-ROM were similar, both costs being orders of magnitude lower than the computational cost of the full order model. 
We also compared the new CF-ROM with state-of-the-art LES-ROMs.
The CF-ROM was as accurate as state-of-the-art LES-ROMs.
However, the CF-ROM was significantly more efficient than these LES-ROMs. 

Although these preliminary results are encouraging, the new CF-ROM framework's full potential still needs to be explored.
Next, we outline several research directions that could be pursued. \\[-0.3cm]

\underline {\it Quadratic Ansatz}: \ 
In the CF-ROM's numerical investigation in Section~\ref{sec:numerical-results}, we used the linear ansatz (i.e., we considered $\tB_{i m n} = 0$ in~\eqref{eqn:c-les-rom-2}) to model the ROM stress tensor $\btau_r$.
Of course, this ensured a simple numerical treatment of the resulting linear least squares problem in Section~\ref{sec:cf-rom}.
We emphasize, however, that from a physical point of view, the quadratic ansatz (i.e., $\tB_{i m n} \neq 0$ in~\eqref{eqn:c-les-rom-2}) is more appropriate.
Indeed, the ROM stress tensor $\btau_r$ is closely connected to the quadratic nonlinearity in the NSE.
Thus, we expect the ROM stress tensor be best approximated by a quadratic term, too, just as in the quadratic ansatz.
Although using the quadratic ansatz will result in a more challenging nonlinear least squares problem, we believe that this will significantly increase the CF-ROM's physical accuracy.
In fact, given the low physical accuracy of the linear ansatz used in the numerical investigation in Section~\ref{sec:numerical-results}, the fact the CF-ROM results were qualitatively comparable to those from state-of-the-art LES-ROMs is impressive. \\[-0.3cm]

\underline {\it CF-ROM's Potential Generality}: \ 
Although the particular form of the CF-ROM will change with the underlying physical system, certain classes of physical systems might yield similar types of CF-ROMs.
For example, for some classes of systems the linear ansatz could suffice, whereas for some other classes of systems the quadratic ansatz could be needed.
Furthermore, for certain classes of systems, one could determine variation intervals for the CF-ROM coefficients. \\[-0.3cm]

\underline {\it CF-ROM's Universality}: \ 
We constructed and tested the CF-ROM in a fluid dynamics setting: 
We used the NSE to build the CF-ROM and the Burgers equation to numerically investigate it.
We emphasize, however, that the CF-ROM framework can be applied to any type of nonlinear PDE that is amenable to reduced order modeling. 
Indeed, the only input needed in the CF-ROM framework is the FOM data. 
Once those are supplied, the CF-ROM proceeds in two steps.
(i) First, the given nonlinear PDE is spatially filtered.
The nonlinearity yields a nonlinear stress tensor (which will generally be different from the stress tensor $\btau_r$ used in this paper).
(ii) In the second step of the CF-ROM construction, the available FOM data is used to compute an approximation for the true stress tensor in the filtered ROM in (i) and an optimization problem is employed to calibrate the CF-ROM accordingly.
We emphasize again that the entire CF-ROM procedure does not use any phenomenological arguments that would restrict it to the particular physical system modeled by the given nonlinear PDE.
This is in stark contrast with, e.g., ROM closure models of eddy viscosity type~\cite{HLB96,osth2014need,wang2012proper}, which cannot be directly applied to other classes of PDEs.
Since the CF-ROM is built upon general principles (i.e., filtering and calibration), we expect it to be successful in the numerical simulation of general mathematical models (e.g., from elasticity or bioengineering). \\[-0.3cm]

\underline {\it Spatial Filtering}: \ 
In Section~\ref{sec:f-rom}, we used the ROM projection~\eqref{eqn:rom-projection} as a ROM spatial filter to construct the F-ROM~\eqref{eqn:f-rom}.
This choice of ROM spatial filter was essential for being able to write the F-ROM decomposition in~\eqref{eqn:les-rom-9} and to explain why most ROM closure models amount to adding extra terms to the standard G-ROM.
(Indeed, the ROM projection yielded identity~\eqref{eqn:les-rom-4b}, which in turn allowed the formal F-ROM decomposition~\eqref{eqn:les-rom-9}.)
We emphasize, however, that other ROM spatial filters could be used in the new CF-ROM framework.
For example, the ROM differential filter (which was successfully used in developing LES-ROMs~\cite{wells2017evolve,xie2017approximate}) could be used as a ROM spatial filter to construct the F-ROM~\eqref{eqn:f-rom}.

\bibliographystyle{plain}
\bibliography{traian}

\end{document}